\documentstyle[amssymb,12pt]{amsart}
\newtheorem{prop}{Proposition}[section]
\newtheorem{defin}{Definition}[section]

\newtheorem{thm}{Theorem}[section]

\theoremstyle{remark}

\begin{document}
\newcommand{\nc}{\newcommand}
\nc{\on}{\operatorname}
\nc{\pa}{\partial}
\nc{\cA}{{\cal A}}
\nc{\cB}{{\cal B}}
\nc{\cC}{{\cal C}}
\nc{\cQ}{{\cal Q}}
\nc{\cE}{{\cal E}}
\nc{\cG}{{\cal G}}
\nc{\cH}{{\cal H}}
\nc{\cX}{{\cal X}}
\nc{\cR}{{\cal R}}
\nc{\cL}{{\cal L}}
\nc{\cK}{{\cal K}}
\nc{\sh}{\on{sh}}
\nc{\Id}{\on{Id}}
\nc{\Diff}{\on{Diff}}
\nc{\ad}{\on{ad}}
\nc{\Der}{\on{Der}}
\nc{\End}{\on{End}}
\nc{\res}{\on{res}}
\nc{\ddiv}{\on{div}}
\nc{\card}{\on{card}}
\nc{\Jac}{\on{Jac}}
\nc{\Imm}{\on{Im}}
\nc{\limm}{\on{lim}}
\nc{\Ad}{\on{Ad}}
\nc{\ev}{\on{ev}}
\nc{\Hol}{\on{Hol}}
\nc{\Det}{\on{Det}}
\nc{\de}{\delta}
\nc{\si}{\sigma}
\nc{\ve}{\varepsilon}
\nc{\al}{\alpha}
\nc{\CC}{{\mathbb C}}
\nc{\ZZ}{{\mathbb Z}}
\nc{\NN}{{\mathbb N}}
\nc{\KK}{{\mathbb K}}
\nc{\PP}{{\mathbb P}}
\nc{\RR}{{\mathbb R}}
\nc{\WW}{{\mathbb W}}
\nc{\zz}{{\mathbf z}}
\nc{\AAA}{{\mathbb A}}
\nc{\cO}{{\cal O}}
\nc{\cS}{{\cal S}}
\nc{\cF}{{\cal F}}
\nc{\la}{{\lambda}}
\nc{\G}{{\mathfrak g}}
\nc{\A}{{\mathfrak a}}
\nc{\HH}{{\mathfrak h}}
\nc{\N}{{\mathfrak n}}
\nc{\B}{{\mathfrak b}}
\nc{\La}{\Lambda}
\nc{\g}{\gamma}
\nc{\eps}{\epsilon}
\nc{\wt}{\widetilde}
\nc{\wh}{\widehat}
\nc{\bn}{\begin{equation}}
\nc{\en}{\end{equation}}
\nc{\SL}{{\mathfrak{sl}}}

\newcommand{\ldar}[1]{\begin{picture}(10,50)(-5,-25)
\put(0,25){\vector(0,-1){50}}
\put(5,0){\mbox{$#1$}}
\end{picture}}

\newcommand{\lrar}[1]{\begin{picture}(50,10)(-25,-5)
\put(-25,0){\vector(1,0){50}}
\put(0,5){\makebox(0,0)[b]{\mbox{$#1$}}}
\end{picture}}

\newcommand{\luar}[1]{\begin{picture}(10,50)(-5,-25)
\put(0,-25){\vector(0,1){50}}
\put(5,0){\mbox{$#1$}}
\end{picture}}

\begin{flushright}
ITEP-TH-42/01\\
hep-th/0110032\\
\end{flushright}

\title[Polynomial Poisson algebras]
{Polynomilal Poisson algebras with regular structure of symplectic leaves}

\author{Alexandre Odesskii}

\address{A.O.:Landau Institute of Theoretical Physics, 2, Kosygina str.,117334 Moscow , Russia}

\address{ D\'ept. de Math\'ematiques, Univ. d'Angers, 2, Bd. Lavoisier, 49045 Angers, France}

\author{Vladimir Rubtsov}

\address{V.R.: D\'ept. de Math\'ematiques, Univ. d'Angers, 2, Bd. Lavoisier, 49045 Angers, France}

\address{ITEP, 25 Bol. Cheremushkinskaia, 117259 Moscou, Russie}

\date{july 2001}

\begin{abstract}

We study polynomial Poisson algebras with some regularity conditions. Linear
(Lie-Berezin-Kirillov) structures on dual spaces of semi-simple Lie algebras,
quadratic Sklyanin elliptic algebras of \cite{FO1},\cite{FO2} as well as polynomial
algebras recently described by Bondal-Dubrovin-Ugaglia (\cite{Bondal},\cite{Ug})
belong to this class. We establish some simple determinantal relations between the
brackets and Casimirs in this algebras. These relations imply in particular that for
Sklyanin elliptic algebras the sum of Casimir degrees coincides with the dimension of
the algebra. We are discussing some interesting examples of these algebras and in
particular we show that some of them arise naturally in Hamiltonian integrable
systems. Among these examples is a new class of two-body integrable systems admitting
an elliptic dependence both on coordinates and momenta.

\end{abstract}

\maketitle

\section{Introduction}

We shall understand under polynomial Poisson structures those ones whose brackets are
polynomial in terms of local coordinates on underlying Poisson manifold. The typical
example of a such structure is the famous Sklyanin algebra.

Remind that a Poisson structure on a manifold $M$(for instance it does not play an
important role is it smooth or algebraic) is given by a bivector antisymmetric tensor
field $\pi\in \Lambda^2(TM)$ defining on the corresponded algebra of functions on $M$
a structure of (infinite dimensional) Lie algebra by means of the Poisson brackets
$$
\{f,g\} = \langle{\pi, df\wedge dg}\rangle.
$$
The Jacobi identity for this brackets is equivalent to an analogue of (classical)
Yang-Baxter equation namely to the "Poisson Master Equation": $[\pi,\pi] = 0$, where
the brackets $[,]:\Lambda^p(TM)\times \Lambda^q(TM)\mapsto \Lambda^{p+q-1}(TM)$ are
the only Lie super-algebra structure on $\Lambda^{.}(TM)$ given by the so-called
Schouten brackets. We refer for all this facts for example to the book
\cite{Vais}.

There is an analogue of the Darboux theorem \cite{Wein} which gives a
local description of any Poisson manifold. Namely,  there are
coordinates $(q_1,\ldots,q_l,p_1,\ldots,p_l,x_1,\ldots,x_k)$ near any
point $m\in M$ such the bivector field $\pi$ reads as
$$
\pi = \sum_{i=1}^{l}{{\partial}\over{\partial
    q_i}}\wedge{{\partial}\over{\partial p_i}} + \sum_{i>j}^{k}f_{ij}(x){{\partial}\over{\partial
    x_i}}\wedge{{\partial}\over{\partial x_j}}
$$
such that $f(m) = 0$.

The case $k=0$ is corresponded to a symplectic structure and the opposite case $l=0$
is called usually totally degenerated Poisson structure. Let us remind the important
notion of the Casimir functions of $\pi$ (or briefly speaking {\it Casimirs}). A
function $F\in Fun(M)$ is a Casimir of the Poisson structure $\pi$ if $\{F,G\} = 0$
for all functions $G\in Fun(M)$. It is clear that if the rank of the structure is
constant in this neighbourhood of $m$ ($m$ is called a regular point) then the
Casimirs in the neighbourhood  are the functions depending only on $x_1,\ldots,x_k$
and Poisson manifold admits a foliation by symplectic leaves, i.e. is a unification of
submanifolds

$$
x_1 = c_1,\ldots,x_k = c_k \,
$$
and $c_i$ are constants

such that $\pi$ is non- degenerate on each of them.
In general the dimension of the leaves is constant only the open dense
subset of regular points and may vary outside.

We will discuss in this paper such polynomial Poisson structures  which have
some regularity properties for  their symplectic leaves. We will postpone
the exact definition of these properties but remark that there are lot
of interesting classes of Poisson structures which acquire them.
The most familiar is of course the famous Lie- Berezin - Kirillov
Poisson structure on the dual space $\G^*$ to a semi-simple Lie algebra $\G$.
This structure is linear in the coordinate functions $x_1,\ldots,x_n$
on  $\G^*$ corresponding to a base $X_1,\ldots,X_n$ of  $\G$ and is
given by the structure constants of $\G$. Reciprocally, any linear
structure
$$
\{x_i,x_j\} = C_{ij}^{k}x_k
$$
arises from a Lie algebra. We should restrict to the semi-simple Lie algebras if we
want to get a regular Poisson Lie-Berezin - Kirillov bracket.

The next wide class of interesting Poisson structures with regular symplectic leaves
is provided by a subclass of {\it quadratic} Poisson algebras introducing by E.
Sklyanin in \cite{sklyanin} and in more general context describing in
\cite{FO1},\cite{FO2} as a quasi-classical limit of associative algebras with
quadratic relations  which are flat deformations of function algebra on $\CC^n$. These
algebras are associated with elliptic curves and have the following curious property:
they have polynomial Casimirs and for all of them the dimension of the algebra is
equal to the sum of the degrees of Casimir polynomial ring generators. One of the aim
of the paper is to give an explanation to this property which follows from a simple
general determinantal relation between alternated products of the coordinate function
Poisson brackets and the minors of the Jacobi matrix of Casimirs with respect to this
coordinates (theorem 3.1 below).

Another motivation to study this type of Poisson polynomial algebras comes from the
world of Hamiltonian integrable systems. We will show that such Poisson algebras arise
naturally as Hamiltonian structures attached to different physically interesting
models: from the trivial Arnold-Euler-Nahm top till the recently described
(\cite{FGNR},\cite{BMMM}) double-elliptic $SU(2)$-model. The polynomial Poisson
algebra which provide a natural Hamiltonian structure to this model is homogeneous and
has the degree 3.

The paper is organized as follows. We start with the general remarks concerning
the polynomial Poisson  algebra. We remind some well (and maybe less)
known definitions, discuss the holomorphic prolongation of the affine
Poisson structures and remind the generalized Sklyanin elliptic
algebras,  following \cite{FO1}.

In the next section we define the algebras with a regular structure of
symplectic leaves and prove the mentioned theorem (Theorem 4.1) which
is just a straightforward calculation. The less evident example of the
theorem 4.1 relation is given by the generalized Sklyanin algebra
$q_5$.

The chapters 3 and 4 are devoted to different examples. We discuss in
ch.3 the examples of amazing {\it non- polynomial} changes of variables which preserve the polynomial character of the Poisson structures associated with
elliptic curves. This changes are going back to "Mirror
transformation" of Calabi-Yau manifolds. The form and sens of this
change in "quantum" setting of \cite{FO1} is an interesting open
question.

In the chapter 4 we discuss the examples of polynomial Poisson
structures associated with K3 surfaces (both smooth or singular). 
The interesting feature of this structures is that they are of degree
3 in coordinate functions and it is still unclear how to "quantizë
them or what are the analogues of generalized Sklyanin algebras to
this case? This analogues (if they are exist ) should be called "Mukai
algebras". At the moment the only known examples are the "quantum"
projective plane of \cite{VdB} or some "non-commutatives" K3 surfaces
(\cite{KimLee}, which are the result of deformation of toroidal
orbifolds 
and don't provide an interesting example of Poisson algebra
quantization).

Finally, in the last section we give the examples of some natural
integrable systems
which are Hamiltonian with respect to the regular Poisson structures
under discussion. In fact, one of the examples, the double elliptic or
"DELL" system which is cubic in the natural Hamilton setting, gave the
initial motivation to study this series of structures. The detailed
study of this interesting example associated with 6d SUSY gauge theory
with matter hypermultiplets will be done in our joint paper with
H.Braden and A. Gorsky. Our paper may be considered as an account to
the polynomial Poisson structures which  will be used for the DELL. This is a
particular explanation of the illustrative character of the paper. We often
had restricted on the level of the evident examples to do the paper
accessible to a wide circle of readers (in the first for physicists).

\section{Poisson algebras on polynomial rings}

\subsection{General remarks}

We will discuss some examples of polynomial Poisson algebras which are defined on an affine part
of some algebraic varieties. Typically, this varieties embedded as (complete)
intersections in (weighted) projective spaces. We should remark that the
considerations of intersection varieties in weighted projective spaces are equivalent
sometimes to intersections in the products of usual ${\CC}P^{n}$ - the fact which is
well -known to string theory physicists who had studied the Calabi -Yau complete
intersections as a compactification scenario to superstring vacua and the mirror
symmetry as a map between the moduli spaces of weighted Calabi - Yau manifolds
(\cite{CDLS},\cite{CGH}).

Let us consider $n-2$ polynomials $Q_i$ in ${\CC}^{n}$ with coordinates
$x_i, i=1,...,n$ . For any polynomial $\lambda \in {\CC}[x_1,...,x_{n}]$ we can define a
bilinear differential operation
$$
\{,\}:{\CC}[x_1,...,x_{n}]\otimes{\CC}[x_1,...,x_{n}] \mapsto
{\CC}[x_1,...,x_{n}]
$$
by the formula
\begin{equation}
\label{bracket} \{f,g\}= \lambda\frac{df\wedge dg\wedge dQ_1\wedge ...\wedge
dQ_{n-2}}{dx_l\wedge dx_2\wedge ...\wedge dx_{n}}, \ f,g \in {\CC}[x_1,...,x_{n}].
\end{equation}

This operation gives a Poisson algebra structure on ${\CC}[x_1,...,x_{n}]$ as a
partial case of more general $n-m$-ary Nambu operation given by an
antysymmetric $n-m$-polyvector field $\eta$ :
$$
\langle \eta,df_1\wedge ...\wedge df_{n - m}\rangle =
\bigl\{f_1,...,f_{n-m}\bigr\},
$$
depending on $m$ polynomial "Casimirs" $Q_1,\ldots,Q_m$ and $\lambda$

such that
\begin{equation}
\label{Nambu}\{f_1,\ldots,f_{n-m}\}= \lambda\frac{df_1\wedge\ldots\wedge df_{n-m}\wedge dQ_1\wedge ...\wedge
dQ_{m}}{dx_l\wedge dx_2\wedge ...\wedge dx_{n}}, \ f_{i}\in {\CC}[x_1,...,x_{n}].
\end{equation}

and

$$
\bigl\{,...,\bigr\}: {\CC}[x_1,...,x_{n}]^{\otimes {n-m}}{\mapsto
{\CC}[x_1,...,x_{n}]}
$$

such that the three properties are valid:

1)antisymmetricity:
$$
\bigl\{f_1,...,f_{n-m}\bigr\} =
(-1)^{\sigma}\bigl\{f_{\sigma(1)},...,f_{\sigma(n-m)}\bigr\}, \sigma
\in Symm_{n-m};
$$

2) coordinate-wise "Leibnitz rule" for any $h\in {\CC}[x_1,...,x_{n}]$:
$$
\bigl\{f_1h,...,f_{n-m}\bigr\} = f_1\bigl\{h,...,f_{n-m}\bigr\} + h\bigl\{f_1,...,f_{n-m}\bigr\};
$$

3) The "Fundamental Identity" (which replaces the Jacobi):
$$
\bigl\{\bigl\{f_1,...,f_{n-m}\bigr\},f_{n-m+1},...,f_{2(n-m)-1}\bigr\}+
$$
$$
\bigl\{f_{n-m},\bigl\{f_1,...,({f_{n-m}})^{\vee}f_{n-m+1}\bigr\},f_{n-m+2},...,f_{2(n-m)-1}\bigr\} +
$$
$$
+
\bigl\{f_{n-m},...,f_{2(n-m)-2},\bigl\{f_1,...,f_{n-m-1},f_{2(n-m)-1}\bigr\}\bigr\}
=
$$
$$
\bigl\{f_1,...,f_{n-m-1},\bigl\{f_{n-m},...,f_{2(n-m)-1}\bigr\}\bigr\}
$$
for any $f_{1},...,f_{2(n-m)-1}\in {\CC}[x_1,...,x_{n}]$.

This structure is a natural generalization of the Poisson structure (which corresponds
to $n-m=2$) was introduced by Y.Nambu in 1973(\cite{nambu}) and was recently extensively
studied by L.Tachtajan \cite{tacht}.

The most natural example of the Nambu - Poisson structure is so-called "canonical"
Nambu-Poisson structure on ${\CC}^m$ with coordinates $x_1,...,x_m$:
$$
\bigl\{f_{1},...,f_{m}\bigr\}= Jac(f_{1},...,f_{m})={{\partial (f_1,...,f_m
)}\over{\partial(x_1,...,x_m)}}.
$$

We should remark also that the formula (\ref{bracket}) for the Poisson
brackets takes place in more general setting when the polynomials
$Q_i$ are replaced, say, by rational functions but the resulting
brackets are still polynomials. More generally, this formula is valid
for the power series rings.

The polynomials $Q_i,i=1,...,n-2$ are Casimir functions for the bracket
(\ref{bracket}) and any Poisson structure in $\CC^{n}$ with $n-2$ generic
Casimirs $Q_i$ are written in this form.

All this facts are not new and very well known probably since Nambu
and Sklyanin papers. They had reappared recently
(\cite{Perel},\cite{Khimsh}) in the different frameworks because of
general interest to classification and quantization problems for
Poisson structures.

The case $n=4$ in (\ref{bracket}) corresponds to the classical (generalized) Sklyanin quadratic
Poisson algebra. The very Sklyanin algebra is associated with the following two quadrics in
$\CC^4$:
\begin{equation}
\label{sklyanin1}
Q_1 = x_1^2 + x_2^2 + x_3^2,
\end{equation}
\begin{equation}
\label{sklyanin2}
Q_2 = x_4^2 + J_1x_1^2 + J_2x_2^2 + J_3x_3^2.
\end{equation}

The Poisson brackets (\ref{bracket}) with $\lambda=1$ between the affine coordinates looks as follows

\begin{equation}
\label{jacobian}
\{x_i,x_j\}=(-1)^{i+j}det\left(\frac{\pa Q_k}{\pa x_l}\right), l\neq
i,j, \, i>j.
\end{equation}

The expression (\ref{bracket}) has an advantage before
(\ref{jacobian}) because it is compatible with the more general
situations when the intersected varieties are embedded in the weighted
projective spaces or in the product of the projective spaces.
We will consider an example of such situation below.

The natural question  arises: to extend the brackets (\ref{bracket}) or
(\ref{jacobian}) from $\CC^{n}$ to the projective space ${\CC}P^{n}$.

We can state the following

\begin{prop}

Let $X_1,...,X_n$ are coordinates on ${\CC}^n$ considering as an
affine part of the corresponding projective space ${\CC}P^n$ with the
homogeneous coordinates $\left(x_0:x_1:\cdots :x_n\right), X_i = \frac{x_i}{x_0}$
then if

$\{X_i,X_j\}$ extends to a holomorphic Poisson structure on
${\CC}P^{n}$ then the maximal degree of the
structure (= the length of monomes in $X_i$)  is $3$ and
\begin{equation}
\label{holomorph}
X_k\{X_i,X_j\}_3 + X_i\{X_j,X_k\}_3 + X_j\{X_k,X_i\}_3 = 0, i\neq j \neq k,
\end{equation}
i.e.  $\{X_i,X_j\}_3 =
X_iY_j -  X_jY_i$, with  $deg Y_i = 2$
\end{prop}

Indeed, under the change of the coordinates we have
$$
\{1/X_i,X_j/X_i\} = -1/{X_i}^3\{X_i,X_j\}
$$
hence $deg\{X_i,X_j\}\leq 3$.

Moreover, we can conclude that
$$
\{X_i,X_j\} = \{X_i,X_j\}_0 + \{X_i,X_j\}_1 + \{X_i,X_j\}_2 + \{X_i,X_j\}_3,
$$
where $deg \{X_i,X_j\}_k = k , k= 0,1,2,3$ .

 The \ref{holomorph} is a corollary of the following identities:
$$
\{X_j/X_i,X_k/X_i\} = 1/{X_i}^2\{X_j,X_k\} - X_k/{X_i}^3\{X_j,X_i\} -
X_j/{X_i}^3\{X_i,X_k\},
$$

The statement is proved.

In general, the homogeneous Poisson algebras which are described by
(\ref{bracket}) or (\ref{jacobian}) have no projective extensions because
they don't satisfy for $n \geq 4$ to the conditions of the proposition.

The Poisson algebras of the type (1), (5) have
the following property: they have only symplectic leaves of "small"
dimensions ("small" means in fact 0 and 2).

The following description of "small" dimensional symplectic leaf structures is a
fairly direct corollary of the leaf definition

\begin{prop}

Let $\{x_i,x_j\} = \langle{\pi,dx_i\wedge dx_j}\rangle = p_{ij}$ be an affine Poisson
structure. Then $\pi$ has only "small dimenisional" symplectic leaves iff $p_{ij}$ is
written in "Plucker form" or  $p_{ij} = \alpha_i\beta_j - \alpha_j\beta_i$ where
$\alpha_i, \beta_j$ are some functions (not necessarily polynomials).
\end{prop}

\subsection{Poisson algebras associated to elliptic curves.}

Another wide class of the polynomial Poisson algebras  arises as a quasi-classical
limit $q_{n,k}(\cE)$ of the associative quadratic algebras $Q_{n,k}(\cE,\eta)$ which
were introduced in a cycle of papers \cite{FO1,FO2}. Here $\cE$ is an elliptic curve
and $n,k$ are integer numbers without common divisors ,such that $1\leq k < n$ while
$\eta$ is a complex number and $Q_{n,k}(\cE,0) = \CC[x_1,...,x_n]$.

Let $\cE = \CC/\Gamma$ be an elliptic curve defined by a lattice $\Gamma =
\ZZ\oplus\tau\ZZ, \tau \in {\CC}, \Im \tau > 0$. The algebra $Q_{n,k}(\cE,\eta)$ has
generators $x_i, i\in {\ZZ}/n{\ZZ}$ subjected to the relations
$$
\sum_{r\in{\ZZ}/n{\ZZ}}{\frac{\theta_{j-i+r(k-1)}(0)}{\theta_{j-i-r}(-\eta)\theta_{kr}(\eta)}}x_{j-r}x_{i+r}
= 0
$$

and have the following properties:

1) $Q_{n,k}(\cE,\eta) = {\CC}\oplus Q_1\oplus Q_2\oplus ...$ such that
$Q_{\alpha}*Q_{\beta}=Q_{\alpha + \beta}$, here $*$ denotes the algebra
multiplication. In other words, the algebras $Q_{n,k}(\cE,\eta)$ are
$\ZZ$ - graded;

2) The Hilbert function of $Q_{n,k}(\cE,\eta)$ is $\sum_{\alpha \geq
  0}\dim Q_{\alpha}t^{\alpha} = \frac{1}{(1 - t)^n}$.

We see that the algebra $Q_{n,k}(\cE,\eta)$ for fixed $\cE$ is a flat
deformation of the polynomial ring $\CC[x_1,...,x_n]$.

Let $q_{n,k}(\cE)$ be the correspondent Poisson algebra. It is shown in \cite{FO1}
that the algebra $q_{n,k}(\cE)$ has $l = \gcd(n,k+1)$ Casimirs. Let us denote them by
$P_\alpha, \alpha \in {\ZZ}/l{\ZZ}$. Their degrees $\deg P_\alpha$ are equal to $n/l$.

We should stress, that the algebras $q_{n,k}(\cE)$ for $n > 4$ have symplectic leaves
of dimension $>2$ and hence these algebras are not given by the formulas
(\ref{bracket}), (\ref{jacobian}) for $n>4$. The examples of the algebras
$q_{3,1}(\cE), q_{4,1}(\cE)$ will be discussed below.

\section{Algebras with regular structures of symplectic leaves}

Let us remind that the minimal codimension of symplectic leaves of a Poisson algebra is called a rank of this algebra.
Now we will describe the class of Poisson algebras on $\CC^n$ satisfying to the following properties of regularity:
\begin{defin}
\label{regular}
A Poisson polynomial algebra $A$ on $\CC^n$ of rank $l$ is called {\it an algebra with a regular structure of
symplectic leaves} if:

1) The center $Z(A)$ is a polynomial ring $\CC[Q_1,\cdots,Q_l]$ of Casimirs $Q_i$;

2) The subvariety $L_{\lambda_1,\cdots,\lambda_l} := \bigcap_i \{Q_i = \lambda_i \}$
is a complete intersection (or the tangent spaces $T_p(\{Q_i = \lambda_i\})$ in
generic point $p\in L_{\lambda_1,\cdots,\lambda_l}$ are intersected transversally);

3) Let $M^*$ be a unification $\bigcup_i \{F_i\}$ of symplectic leaves of dimension
$\dim F_i < n-l$
   then $\dim M^* \leq n-2$;

4) Let $L^*$ be a unification of singularities of $L_{\lambda_1,\cdots,\lambda_l}$ by $\lambda_1,\cdots,\lambda_l$,
   then $\dim L^* \leq n-2$.
\end{defin}

We want to remark that the elliptic algebras  $q_{n,k}(\cE)$ are satisfied to the definition as follows
 from the description of their symplectic leaves in \cite{FO2}.
 For all of them as well as for the others algebras with the property \ref{regular} the following theorem holds

\begin{thm}
\label{hypotesa}
$$
\underbrace{\pi\wedge\pi\wedge\cdots\wedge \pi}_{\mbox{(n-l)/2}} =
\lambda (dQ_1\wedge\cdots\wedge dQ_l)^*,\, \lambda \in {\CC^*},
$$
where $l$ is the dimension of the Poisson center (the number of Casimirs) and
$Q_1,\cdots,Q_l$ are the Casimirs and $*$ means duality between $l$-forms and
$l$-polyvectors established by the standard choice of a volume form
$dx_1\wedge\cdots\wedge dx_n$.
\end{thm}

In coordinates the expression (\ref{hypotesa}) is re-written as
$$
\lambda\det
\left(
\begin{array}{ccc}
\label{hypotesacoord}
{\frac{\pa Q_1}{\pa x_{i_{n-l+1}}}}&\ldots &{\frac{\pa Q_1}{\pa x_{i_{n}}}}\\
\vdots& \ddots & \vdots\\
{\frac{\pa Q_l}{\pa x_{i_{n-l+1}}}}&\ldots &{\frac{\pa Q_l}{\pa x_{i_{n}}}}
\end{array}
\right) = Alt_{(i_1,\cdots,i_n)}\left(\{x_{i_1},x_{i_2}\}\cdots\{x_{i_{n-l-1}},x_{i_{n-l}}\}\right),
$$
where $(i_1,\cdots i_n)$ is an even permutation of $(1,\cdots,n),
\lambda \in {\CC^*}-$non-zero constant.

{\bf Proof of \ref{regular}} For any $n-l$ polynomial  functions $f_1,\cdots,f_{n-l}$  the formula \ref{regular} may be
 re-written in the following way
\begin{equation}
\label{lambda}
w^{(n-l)/2}(f_1,\cdots,f_{n-l}) =
\lambda\frac{df_1\wedge\cdots\wedge df_{n-l}\wedge dQ_1\wedge\cdots\wedge dQ_l}{dx_1\wedge\cdots\wedge dx_n},\
 \lambda\in {\CC^*}.
\end{equation}

Let $\Omega_1$(resp. $\Omega_2$) be left(resp. right) part of (\ref{lambda}).

It is clear that $\imath_{Q_{\alpha}}\Omega_i = 0$ for any $\alpha, i=1,2$, and hence writing the
 tensors $\Omega_i$ in coordinates
$$
\Omega_i = \sum_{(\alpha_1,\cdots\,\alpha_{n-l})} w_{\alpha_1,\cdots,\alpha_l}^{i}\pa_{\alpha_1}\wedge\cdots,\wedge \pa_{\alpha_{n-l}}
$$
we obtain a system of linear equations to the coefficients $w_{\alpha_1,\cdots,\alpha_l}^{i}$:

$$
\imath_{Q_{\alpha}}\Omega_i=
\sum_{\alpha_1}w_{\alpha_1\cdots\alpha_{n-l}}^{i}Q_{\alpha x_{\alpha_1}}^{'}\pa_{\alpha_2}\wedge\cdots\wedge\pa_{\alpha_{n-l}} =0.
$$

We obtain that both sides of the (\ref{lambda}) are proportional one to
another (where $\lambda$ is a non-zero function)
from this system.
This function $\lambda$ is rational and hence there are two
polynomials without common divisors $p_1,p_2$  such that
$\lambda = \frac{p_1}{p_2}$.
But the first polynomial $p_1$ has the set of zeroes lying in $M^*$ hence $codim M^* = 1$.
Taking in to the account 3) of \ref{hypotesa} we obtain $p_1 = const$.

The same arguments show that the zeroes of $p_2$ contain in (or coincide with) $L^*$ and $p_2 = const_1$.
 Hence the result.

{\bf Remark 1} The statement of \ref{hypotesa} is still valid to {\it any} Poisson algebra of rank $l$ with $l$
independent central elements $Q_{\alpha}$ if we admit as the factor $\lambda$ an arbitrary rational function.

{\bf Remark 2} Under the  conditions of the theorem if $A$ is a quadratic Poisson
algebra our formula gives the following relation between the degrees of Casimirs
and the dimension of $A$:
$$
\sum_{\alpha}\deg Q_{\alpha} = dim A = n
$$
Moreover , for all known
elliptic algebras which are graded deformations of ${\ZZ}^n$-graded polynomial ring
with finite number of generators, the sum of generator degrees is equal to the sum of
degrees of Casimirs, For example, $n = \sum_{i\in {\ZZ/n\ZZ}}\deg x_i = \sum_{j\in
{\ZZ/l\ZZ}}\deg Q_j$.

{\bf Remark 3} In the case when the condition 3) in \ref{hypotesa} is failed
the following definition is useful:
\begin{defin}
An element $Q \in A$ is called {\it quasi-central} or quasi-Casimir if the
bracket $\{Q,f\}$ is divided by $Q$ for any $f\in A$.
\end{defin}
It is clear that if $Q_1$ and $Q_2$ are quasi-central elements of $A$ then their
product $Q_1Q_2$ is also quasi-central. Let $codim M^* = 1$ and $\{M_i\}$ are
irreducible components of $M^*$ of dimension $n-1$. If $M_i$ is given by the equations
$Q_i = 0$, then it is clear that $\{Q_i\}$ are quasi-Casimirs. Conversely, let
undecomposable $Q$ be a quasi-central but non-central element, then the hyper-surface
$Q = 0$ contains in $M^*$.

The last curious and possibly useful observation is that for any quasi-central $Q$ the bracket
$\{log Q,.\}$ is an (exterior) derivation of $A$.

To illustrate the theorem we show that it is true for the algebra $q_{5}(\cE)$ from
\cite{FO1}.

{\bf Example} We have the polynomial ring with 5 generators $x_i, i\in{\ZZ/5\ZZ}$
enabled with the following Poisson bracket:
$$
\{x_{i+1},x_{i+4}\} = -1/5(k^2 + 3/k^3)x_{i+1}x_{i+4} +
2x_{i+2}x_{i+3} + kx_{i}^2;
$$
$$
\{x_{i+2},x_{i+3}\} = 1/5(3k^2 - 1/k^3)x_{i+2}x_{i+3} +
2/k x_{i+1}x_{i+4} - 1/k^2 x_{i}^2.
$$

Here $i\in \ZZ/5\ZZ$ and $k\in \CC$ is a parameter of the curve $\cE_{\tau}
= \CC/{\Gamma}$, i.e. some function of $\tau$.

The center $Z(q_5(\cE))$ is generated by the polynomial
$$
P = -1/k\left(x_0^5 + x_1^5 + x_2^5 + x_3^5 + x_4^5\right) +
$$
$$
(1/k^5 - 3)\left(x_0^3x_1x_4 + x_1^3x_0x_2 + x_2^3x_1x_3 + x_3^3x_2x_4 +x_2^3x_0x_3\right) +
$$
$$
+ \left(k^3 + 3/k^3)(x_0^3x_2x_3 + x_1^3x_3x_4 +
x_2^3x_0x_4 + x_3^3x_1x_0 + x_4^3x_1x_2\right) -
$$
$$
- (2k + 1/k^4)\left(x_0x_1^2x_4^2 + x_1x_2^2x_0^2 + x_2x_0^2x_4^2
  +x_3x_1^2x_0^2 + x_4x_1^2x_2^2\right) +
$$
$$
+(k^2 - 2/k^3)\left(x_0x_2^2x_3^2 + x_1x_3^2x_4^2 + x_2x_0^2x_4^2 +
  x_3x_1^2x_0^2 + x_4x_1^2x_2^2\right) +
$$
$$
+ (k^4 + 16/k -1/k^6)x_0x_1x_2x_3x_4.
$$

It is easy to check that for any $i\in \ZZ/5\ZZ$
$$
\{x_{i+1},x_{i+2}\}\{x_{i+3},x_{i+4}\} + \{x_{i+3},x_{i+1}\} \{x_{i+2},x_{i+4}\} +
\{x_{i+2},x_{i+3}\}\{x_{i+1},x_{i+4}\} = 1/5\frac{\pa P}{\pa x_i}.
$$

\section{Examples of regular leave algebras,\ $n=3$}

\subsection{Elliptic algebras}
Let
\begin{equation}
\label{torus}
P(x_1,x_2,x_3) = 1/3 \left(x_1^3 + x_2^3 + x_3^3  \right) + kx_1x_2x_3,
\end{equation}
then
$$
\begin{array}{rcl}
\{x_1,x_2\} = kx_1x_2 + x_3^2\\
\{x_2,x_3\} = kx_2x_3 + x_1^2 \\
\{x_3,x_1\} = kx_3x_1 + x_2^2.
\end{array}
$$
The quantum counterpart of this Poisson structure is the algebra
${\cQ}_3({\cE},\eta)$, where ${\cE}\subset {\CC}P^2$ is an
elliptic curve given by $P(x_1,x_2,x_3) = 0$.

\subsection{"Mirror transformation"}

The interesting feature of this algebra is that their
polynomial character is preserved even after the following changes
of  variables:

a)let
\begin{equation}
\label{change} y_1 = x_1, y_2 = x_2x_3^{-1/2}, y_3 = x_3^{3/2}.
\end{equation}

The polynomial $P$ in the coordinates $(y_1,y_2,y_3)$ has the form

\begin{equation}
\label{mirror torus} P^{\vee}(y_1,y_2,y_3)=1/3 \left(y_1^3 + y_2^3y_3 + y_3^3  \right)
+ ky_1y_2y_3
\end{equation}

 and the Poisson bracket is
also polynomial (which is not evident at all!) and has the same form:

$\{y_i,y_j\} = \frac{\pa P^{\vee}}{\pa y_k}$, where $(i,j,k) = (1,2,3)$.

If we put $\deg y_1 = 2,\deg y_2 = 1, \deg y_3 = 3$ then the polynomial $P^{\vee}$ is
also homogeneous in $(y_1,y_2,y_3)$ and defines an elliptic curve $P^{\vee} = 0$ in
the weighted projective space ${\WW}P_{2,1,3}$.

b)now let $z_1 = x_1^{-3/4}x_2^{3/2}, z_2 =  x_1^{1/4}x_2^{-1/2}x_3, z_3 = x_1^{3/2}$.

The polynomial $P$ in the coordinates $(z_1,z_2,z_3)$ has the form $P(z_1,z_2,z_3)
=1/3 \left(z_3^2 + z_1^2z_3 + z_1z_2^3  \right) + kz_1z_2z_3$ and the Poisson bracket
is also polynomial (which is not evident at all!) and has the same form: $\{z_i,z_j\}
= \frac{\pa P}{\pa z_k}$, where $(i,j,k) = (1,2,3)$.

If we put $\deg z_1 = 1,\deg z_2 = 1, \deg z_3 = 2$ then the polynomial $P$ is also
homogeneous in $(z_1,z_2,z_3)$ and defines an elliptic curve $P = 0$ in the weighted
projective space ${\WW}P_{1,1,2}$.

The origins of the strange non-polynomial change of variables (\ref{change}) lie in
the construction of "mirror" dual Calabi - Yau manifolds (\cite{CDLS}) and the torus
(\ref{torus}) has (\ref{mirror torus}) as a "mirror dual". Of course, the mirror map
is trivial for 1-dimensional Calabi - Yau manifolds. Curiously, mapping (\ref{change})
being a Poisson map if we complete the polynomial ring in a proper way and allow the
non - polynomial functions gives rise to a new "relation" on quantum level: the
quantum elliptic algebra $Q_{3}(\cE^{\vee})$ corresponded to (\ref{mirror torus}) has
complex structure $(\tau + 1)/3$  when (\ref{torus}) has $\tau$. Hence, these two
algebras are different. The "quantum" analogue of the mapping (\ref{change}) is still
obscure and needs further studies.

\subsection{Markov polynomials: "rational degeneration"}

We consider as a Casimir now an example of non-homogeneous polynomial ("Markov
polynomial of degree 3") $P(x_1,x_2,x_3) = x_1^2 + x_2^2 + x_3^2 + 3x_1x_2x_3$. The
corresponding Poisson algebra $\{x_i,x_j\} = \frac{\pa P}{\pa x_k}, (i,j,k) = (1,2,3)$
appeared recently as a Hamiltonian structure on the space of Stokes matrices in the
theory of isomonodromic deformations (see \cite{Du},\cite{Ug},\cite{Boalch}) as well
as a quasi-classical limit in quantization of Teichmuller space by \cite{FockCh}. We
should remark that the $n>3$ generalization of this example associated with "higher"
degree invariants studied in full generality by A. Bondal (\cite{Bondal}) have a
regular structure of symplectic leaves.

For example, the next Bondal-Dubrovin-Ugaglia  (\cite{Bondal}, \cite{Ug}) affine
Poisson algebra in $\CC^6$ with two Casimirs of degree 4 and 2 which has a
4-dimensional symplectic leaf given by their intersection satisfies to \ref{regular}:

Let $\CC^6$ is identified with the six-dimensional space of Stokes matrices of the
form
$$
\left(
\begin{array}{cccc}
1& p& q& r\\
0& 1& x& y\\
0& 0& 1& z\\
0& 0& 0& 1
\end{array}
\right)
$$

with Poisson algebra linear-quadratic relations (see \cite{Bondal}).

If we choose two Casimirs
\begin{eqnarray}
& P_1 = p^2 + q^2 + r^2 + x^2 + y^2 + z^2 - pqx - pry -qrz - xyz + prxz\\
& P_2 = pz + xr -qy
\end{eqnarray}

then it is straightforward to check for example the following relation
$$
\{x,y\}\{p,z\} + \{y,z\}\{p,x\} + \{z,x\}\{p,y\} = \det \left(
\begin{array}{cc}
\frac{\pa P_1}{\pa q}& \frac{\pa P_1}{\pa r}\\
\frac{\pa P_2}{\pa q}& \frac{\pa P_2}{\pa r}
\end{array}
\right) =
$$

$$
\det \left(
\begin{array}{cc}
-y& x\\
2q -px -rz& 2r -py - qz + pxz
\end{array}
\right).
$$

\subsection{Polynomial extension of Askey - Wilson algebra}

A classical (Poisson) Askey-Wilson (\cite {KoZhe}) algebra is
described in our terms by the following Casimir polynomial
\begin{equation}
\label{AW}
P(x,y,z) = z^2 - F(x,y), F(x,y) = a x^2 y^2 + g(x,y),
\end{equation}
where $g(x,y) = a_1 x^2 y + a_2 x y^2 + a_3 x^2 + a_4 y^2 +  a_5 xy +
a_6 x + a_7 y$.

The terminology is explained by the fact that the "quantum" analogue
of the algebra admits a natural representations by Askey-Wilson polynomials.

It is easy to obtain a polynomial extension of the Askey-Wilson
algebra with the same rule (\ref {AW})  but with $F(x,y)$ of arbitrary
degree with the only constraint that $P(x,y,z)$ is a Casimir.

For example, if we admit in $F(x,y)$ the maximal degree terms
$x^4,y^4$ we have the extension of Askey-Wilson algebra which is
equivalent to the standard Sklyanin algebra (\ref{sklyanin1}) (\ref{sklyanin2}) (see \cite {KoZhe}).

We assume the locus $P(x,y,z) = 0$ for (\ref{AW}) to be a curve in the
weighted projective space ${\WW}P_{1,1,2}$ with the variables of $deg\ x =
deg\ y = 1, deg\ z =2$. Then the relation of the {\bf Remark 2} in the
theorem 3.1 is still holds.

\section{K3 - surfaces}

\subsection{Cone of a canonical curve}
Let ${\cC}$ be a curve of genus $g > 1$ which is not a hyper-elliptic. Let $K_{\cC}$
be the canonical bundle of ${\cC}$. We consider the canonical embedding $\xi : {\cC}
\to {\CC}P^{g-1}$ and the cone ${\KK}_{\cC}\subset {\CC}^{g}$  of $\xi({\cC})$. Let
${\tilde{\KK}_{\cC}}$ be a covering of ${\KK}_{\cC}$ which correspond to the covering
${\cH}\mapsto {\cE}$ where ${\cH}$ is the upper half-plane with the coordinate $\tau$
and
$$
{\tilde{\KK}_{\cC}} = {\cH}\times{\CC}^*.
$$

We will denote by $d\tau$ the corresponding coordinate on ${\CC}^*$. It is well-known
that ${\KK}_{\cC} ={\tilde{\KK}_{\cC}}/{\Gamma}$ where ${\Gamma}\subset{\SL}_2({\RR})$
is a discrete subgroup which acts naturally on the coordinates $(\tau, d\tau)$. Let us
describe a Poisson structure on ${\tilde{\KK}_{\cC}}$ by the relation $\{\tau,d\tau\}
= (d\tau)^2$. It is easy to check that this structure is ${\SL}_2({\RR})$-invariant
and symplectic. Hence the cone ${\KK}_{\cC}$ is a homogeneous symplectic. For $g=3$
the cone ${\KK}_{\cC}\subset {\CC}^{3}$ is given by the equation $p(x_1,x_2,x_3) = 0$,
where $p(x_1,x_2,x_3)$ is a homogeneous polynomial of the degree 4. We supply
${\CC}^3$ with the polynomial Poisson structure as above: $\{x_i,x_j\} =\frac{\pa
p}{\pa x_k}, (i,j,k) = (1,2,3)$.

\begin{prop}
${\KK}_{\cC}\subset {\CC}^{3}$ is a symplectic leaf of this polynomial structure and
the restriction of the structure to ${\KK}_{\cC}$ has the form
$$
\{\tau,d\tau\} = \lambda (d\tau)^2, \ \lambda\in {\CC}^*.
$$
\end{prop}
We have analogues of this description for the cases $g=4,5$.
\subsection{General discussion.}
Let us compactify the precedent situation in considering the one -point
compactification of the cone or the completion ${\PP}(K_{\cC}\oplus {\cO}_{\cC})$. We
obtain a compact two-dimensional complex manifold with trivial anti-canonical bundle
or a holomorphic embedding of the curve $C$ in a $K3$-surface. We will describe the
restriction of the Poisson structures on the $K3$-surface as a symplectic leave. It is
an interesting example of non-quadratic affine brackets.

All $K3$ have symplectic holomorphic form due to S. Mukai (\cite{mukai}) and all
examples of the projective models of $K3$ has realization of the Mukai symplectic
brackets as the restrictions of polynomial Poisson structures on $\CC^n, n=3,4,5$
arising on the complete intersections of Casimirs which realize these models. These
brackets have the degree 3. It is worst to remark that in the accordance with the
proposition 2.1 the brackets can not be extended to holomorphic Poisson structures on
${\CC}P^n, n=3,4,5$. We will give a direct verification that a  non-trivial
obstruction exists in this case.

It is well-known that the K3-surfaces are the only two-dimensional compact complex
manifolds with zero 1-st Betti number and a holomorphic symplectic structure. We will
give a simple arguments related the polynomial Poisson brackets with Mukai description
of the holomorphic symplectic structure on algebraic K3 surfaces.

The construction of Mukai (\cite{mukai}) is based on the following simple facts: all
algebraic K3 have a projective embedding, for any ${\CC}P^{n}$ with a system of
homogeneous coordinates $(X_0:X_1:...:X_n)$ and the standard $n$-form
$$
\Omega = \sum_{i=0}^{n}(-1)^{i}X_i dX_0\wedge...\wedge{\hat{dX_i}}\wedge...\wedge
dX_n,
$$
which is used to define the Mukai holomorphic two-form by the residue form
$Res_{\cS}\left(\frac{\Omega}{P}\right)$ along the K3 -surface ${\cS}$ where the
projective embedding of ${\cS}$ given by zero loci $P$, which is either quartic
($n=3$) or the transversal intersection of a quadric and a cubic ($n=4$) or the
transversal intersection of three quadrics for $n=5$.

We can check straightforwardly that the polynomial Poisson structure on ${\CC}^n$ for
$n =3,4,5$ given above being restricted to the surface $\cS$ coincides with the Mukai
structure being written in the affine coordinates.

\subsection{Example 1: Fermat quartic}
Let $P_4(X_0,X_1,X_2,X_3) = 0, deg P_4=4$ be a quartic K3 surface in ${\CC}P^3$. Taking an
open domain $U_0 =\{X_0 \neq 0\}$ we have that the form $\alpha =
\frac{\Omega}{X_0^4}$ which is written in the affine part as $\alpha_0 = dx_1\wedge
dx_2\wedge dx_3, x_i = \frac{X_i}{X_0}$ is a holomorphic 3-form on ${\cS}\setminus
U_0$ and has simple poles along the intersection ${\cS}\bigcap U_0$. Hence in the
affine part $\frac{\Omega}{P} = \frac{\alpha_0}{P(1,x_1,x_2,x_3)}$ and the residue
$\omega = Res_{\cS}\frac{\Omega}{P_4}$ is given (for example) by
\begin{equation}
\omega = \frac{dx_1\wedge dx_2}{\frac{\pa P_4}{\pa x_3}}
\end{equation}
and as we have seen  in the affine chart $U_0$ the bracket is given by $\{x_1,x_2\} =
-{\frac{\pa P_4}{\pa x_3}}$.

\begin{prop}
The polynomial structures on $\CC^n, n = 3,4,5$ given by the polynomials $P_4$ for
$n=3$, and by complete intersections of a quadric and cubic $P_2,P_3$ for $n=4$ or by
an intersection of three quadrics $P_2,Q_2,R_2$ for $n=5$ have no holomorphic
extensions to the projective spaces ${\CC}P^n, n = 3,4,5$
\end{prop}

To prove it, for $n=3$ by the proposition 2.1 we can verify that
$$
X_3\{X_1,X_2\}_3 + X_1\{X_2,X_3\}_3 + X_2\{X_3,X_1\}_3 =
-(X_3\frac{\pa P_4}{\pa x_3} + X_1\frac{\pa P_4}{\pa x_1} +
X_2\frac{\pa P_4}{\pa x_2}) = C\neq 0,
$$

Strictly speaking we had checked only that the
structure on $K3$ is not extended by the {\it given} polynomial
formulas. In fact it is easy to
show that the brackets has no {\it any} extension from the
surface to the whole ${\CC}P^3$. Indeed, if such extension exists it
should be polynomial and coincide on $K3$ with our brackets and hence their difference
equals to $0$ on the surface and so it should have in affine
coordinates the degree at least 4, because it should have the polynomial
$P_4$ as a divisor.

\subsection{Example2 : "singular" $K3$ in a product of projective spaces.}
Another interesting example of the formulas to a Poisson structure on $K3$ is based on
the so called "splitting principle " which goes back to the construction of Calabi-Yau
varieties and their mirror dual in \cite{CDLS}. Roughly speaking, the "splitting
principle " permits us to consider the $K3$ surface (possibly singular) embedded as a
hyper-surface in a weighted projective space like a part of two-dimensional variety in
a product of usual (or more generally also weighted) projective spaces. We are able to
apply the coordinate formulas in this situation and we are sure that more rigorous and
conceptual approach is associated with the proper generalization of the residue theory
to toric varieties (see \cite{CCD}).

Let us consider, following to \cite{LS} a hyper-surface
$$
P = y_1(y_1^3 + y_3^6 + y_5^3) - y_2(y_2^3 + y_4^6 - y_5^3) = 0
$$
in ${\WW}P_{1,1,2,2,2}[8]$ which becomes a singular $K3$ in ${\CC}P^3$
\begin{equation}
\label{singular}
P = x_1(x_1^3 + x_3^3 + x_4^3) - x_2(x_2^3 + x_3^3 - x_4^3) = 0
\end{equation}
after the re-definitions $ y_4 = y_3$ and $x_1 = y_1,x_2 = y_2, x_3 = y_3^2,x_4= y_5
$.
This surface in the patch $x_1\neq 0$ is given by the affine equation
\begin{equation}
\label{affsingular}
P = 1 + X_3^3 + X_4^3 - X_2^4 - X_2X_3^3 + X_2X_4^3 = 0
\end{equation}
and in the accordance with the prescription
$$
\{X_2,X_3\} = -\frac{\pa P}{\pa X_4} = -3X_4^2(X_2 + 1).
$$
We can re-write this surface as an intersection

\begin{eqnarray}
& P_1 = z_1x_1 + z_2x_2\\
& P_2 = z_1\left(x_2^3 + x_3^3 - x_4^3\right) + z_2\left(x_1^3 + x_3^3
  + x_4^3\right)
\end{eqnarray}

in the product
\begin{equation}
\begin{array}{c}
 {\CC}P_1\\
 {\CC}P_3
\end{array}
\left[\begin{array}{cc}
1&1\\3&3
\end{array}\right]
\end{equation}

and in the domain $U_{11}\ = \left(\{z_1\neq
  0\}\right)\bigcap\left(\{x_1\neq 0\}\right)$ with the affine
  coordinates $Z ={\frac{z_2}{z_1}}, X_i ={\frac{x_i}{x_1}}, i = 2,3,4$
we have , for example for the intersection
\begin{eqnarray}
& P_1 = 1 + ZX_2 = 0\\
& P_2 = X_2^3 + X_3^3 - X_4^3 + Z\left( 1 + X_3^3 + X_4^3\right) = 0
\end{eqnarray}
$$
\{X_2,X_3\} = \frac{dX_2\wedge dX_3\wedge dP_1\wedge dP_2
}{dX_2\wedge dX_3\wedge dX_4\wedge dZ} = - \frac{dX_2\wedge
dX_3}{\Jac\left(\frac{P_1,P_2}{X_4,Z}\right)} = - 3X_4^2(X_2 + 1).
$$
We can see that this brackets are not homogeneous ( and it is not amazing - we takes
deal with the {\it weighted homogeneous} coordinates - but they have the degree 3 like
it should be for the "usual" $K3$ surfaces.

{\bf Remark} Like it was observed in subsection 4.2 the general definition of the
"splitting principle " leads also to {\it rational}  non -polynomial mapping between
products of  curves,surfaces etc.in products of (weighted) projective spaces. So it is
not clear at the moment which relation it has with the algebras $Q_{n,k}(\cE)$ and
their tensor products. The understanding of a proper sense of the "splitting principle
on "quantum" level may shed some light on the question of a "quantization" of $K3$
surfaces and fibrations.

\section{Integrable systems.}

In this section we collect some examples when the discussed polynomial Poisson
structures naturally appears in Hamiltonian systems. The notion of
Liouville integrability is ambiguous for algebraic varieties and lies out of the scope
of our interests in the paper (see for the discussion \cite{Markush}).We would like only to show that
the naive construction of the affine Poisson brackets sometimes may be very useful and
can clarify the properties of the initial Hamiltonian systems. As a general reference
to the subject we propose the reviews of \cite{AvM}.

\subsection{Quadratic algebras,\ $n=4$}

We will consider a couple of generic quadratic forms in $\CC^4$.
Let
$$
\begin{array}{rcl}
p_1 = 1/2\left(x_1^2 + x_3^2\right) + kx_2x_4\\
p_2 = 1/2\left(x_2^2 + x_4^2\right) + kx_1x_3,
\end{array}
$$

then the brackets are read
$$
\begin{array}{rcl}
\{x_i,x_{i+1}\} = k^2x_{i}x_{i+1} - x_{i+2} x_{i+3}\\
\{x_i,x_{i+2}\} = k\left(x_{i+3}^2 -  x_{i+1}^2\right),
\end{array}
$$

where $i\equiv 0,1,2,3 (mod\ 4).$

The quantization of this relations is the elliptic algebra ${\cQ}_4({\cE},\eta)$,
where  ${\cE}\subset {\CC}P^3$ given by the relations $p_1= 0;p_2 = 0.$

To relate a natural integrable system with this Poisson structure we
need to re-write the couple in the "standard" Sklyanin form
(\ref{sklyanin1}),(\ref{sklyanin2}) and to take the coordinate $x_4$
as a Hamiltonian. Then we obtain from the Sklyanin algebra relations
the following Hamiltonian system:
$$
\begin{array}{rcl}
{\dot x_1} = \{x_1, H\} = (J_2 - J_3)x_2x_3\\
{\dot x_2} = \{x_2, H\} = (J_1 - J_3)x_1x_3\\
{\dot x_3} = \{x_3, H\} = (J_1 - J_2)x_1x_2.
\end{array}
$$
This is so-called "top-like "(or elliptic rotator) representation of the classical
two-particle elliptic S. Ruijsenaars model. This  observation (on the quantum level is
due to I. Krichever and A. Zabrodin (\cite{kriza})) will be discussed in our
subsequent paper with H.Braden and A. Gorsky (\cite{BGOR}).

{\bf Remark} The "usual" Euler - Nahm top is obtained in this scheme as a system on a
co-adjoint orbit of $SO(3,\CC) = SL_2(\CC)$ given by the Casimir $Q_1 = 1/2(x_1^2 +
x_2^2 + x_3^2) = c_1$ with the second quadric $Q_2 = 1/2(J_1 x_1^2 +  J_2x_2^2 +
J_3x_3^2)$ taking as a Hamiltonian. Namely,
$$
\begin{array}{rcl}
{\dot x_1} = \{x_1, H\} = J_2x_2\{x_1,x_2\} +J_3x_3\{x_1,x_3\} = (J_3
-J_2)x_2x_3\\
{\dot x_2} = \{x_2, H\} = J_1x_1\{x_2,x_1\} +J_3x_3\{x_2,x_3\} = (J_1
-J_3)x_1x_3 \\
{\dot x_3} = \{x_3, H\} = J_1x_1\{x_3,x_1\} +J_2x_2\{x_3,x_2\} = (J_2
-J_1)x_1x_2,
\end{array}
$$
where we have used the linear regular Poisson structure associated
with the curve $Q_1 = 1/2(x_1^2 +  x_2^2 + x_3^2) = c_1$ in $\CC^3 :
\{x_i,x_j\} = \epsilon_{ijk}{\frac{\pa Q_1}{\pa x_k}} =
\epsilon_{ijk}x_k.$

\subsection{Double elliptic ( "DELL") system}

There is a deep relation between integrable systems and $N=2$ SUSY gauge theories
which goes back to Witten and Seiberg (\cite{SW}). The classically known
multi-particle integrable Hamiltonian models (both - continues - Calogero-Mother and
differences - Rujisenaars systems)get some new insights coming from the physical
background. The description of 6$d$ gauge theory with the adjoint matter fields had
motivated an appearence of double elliptic (DELL) system unifying the Calogero
-Ruijsenaars family and their {\it duals } and admitting an elliptic dependence both
in coordinates and in momenta in it Hamiltonian (\cite{FGNR}),(\cite{BMMM}). We will
show that there is an example of $K3-$like regular Poisson structure which provides a
natural Hamiltonian description of the $SU(2)$ DELL system.

 The Hamiltonian of the two-particle DELL in center masses frame ($SU(2)$-case)
in the form of (\cite{BMMM}) is given by the function
\begin{equation}
\label{dell}
H(p,q)=\alpha(q|k)cn(p \beta(q|k,\tilde{k})| \frac{\tilde{k}
\alpha(q|k)}{ \beta(q|k) })
\end{equation}
where $\alpha(q|k)=\sqrt{1+ \frac{g^2}{sn^2(q|k)}}$ and $\beta(q|k,\tilde{k})=\sqrt{1+
\frac{g^2\tilde{k}}{sn^2(q|k)}}$, coincides with $x_5$.

We choose the following system of four quadrics in ${\CC}^6$ which provides the phase
space for two-body double elliptic system
$$
\begin{array}{rcl}
x_1^2 - x_2^2 =1 \\
x_1^2 - x_3^2 =k^2 \\
-g^2x_1^2 + x_4^2 - x_5^2=1 \\
-g^2  x_1^2 +  x_4^2 +  \tilde {k}^{-2}x_6^2=\tilde {k}^{-2}
\end{array}
$$
The first pair of the equations yields the "affinization" of projective embedding of
the elliptic curve into ${\CC}P^{3}$ and the second pair provides the elliptic curve
which locally is fibered over the first elliptic curve. If the coupling constant $g$
vanishes the system is just  two copy of elliptic curves embedded in ${\CC}P^{3}\times
{\CC}P^{3}$. Let us emphasize that the coupling constant amounts to the additional non
- commutativity between the coordinates compared to the standard non-commutativity of
coordinates and momenta.

The relevant Poisson brackets for this particular system of quadrics reads
\begin{eqnarray}
\label{quadrics}
&\{x_1,x_2\}= \{x_1,x_3\}=\{x_2,x_3\}=0 \\
&\{x_5,x_1\}= -x_2 x_3 x_4 x_6 \\
& \{x_5,x_2\}= -x_1 x_3 x_4 x_6 \\
& \{x_5,x_3\}= -x_1 x_2 x_4 x_6 \\
& \{x_5,x_4\}= -g^2 x_1 x_2 x_3 x_6 \\
& \{x_5,x_6\}=0
\end{eqnarray}

We should remark that the Poisson structure is singular and can't be extended up to a
holomorphic structure on the whole ${\CC}P^{6}$ because of (\ref{holomorph}).

Due to the commutation relations $x_6$ is constant  which has to be related to the
constant energy to provide the consistency of the system.

The nontrivial commutation relations between coordinates on the distinct tori
correspond to the standard phase space Poisson brackets while the nontrivial bracket
$\{x_5,x_4\}$ means the additional non-commutativity of the momentum space . Let us
note that the triple $x_1, x_2, x_3$ can be considered in the elliptic
parametrization
$$
\begin{array}{rcl}
x_1=\frac{1}{sn(q|k)}\\
x_2=\frac{cn(q|k)}{sn(q|k)}\\
x_3=\frac{dn(q|k)}{sn(q|k)}
\end{array}
$$
by Jacobi sine, cosine and dn functions on the "coordinate" elliptic curve (the torus
with local coordinate $q$).

\begin{thm}
The Hamilton system with the DELL Hamiltonian (\ref{dell}) is equivalent to the
following Hamiltonian system with respect to (\ref{quadrics})and with the Hamiltonian
$x_5$ :
$$
\begin{array}{rcl}
 {\dot x_1} = x_2x_3x_4x_6\\
 {\dot x_2} = x_1x_3x_4x_6\\
 {\dot x_3} = x_1x_2x_4x_6\\
 {\dot x_4} = g^2x_1x_2x_3x_6\\
 {\dot x_6} = 0.
\end{array}
$$
\end{thm}

For the proof see (\cite{BGOR}).

This form of DELL system manifests its algebraic nature and immediately provides its
explicit integration by hyper-elliptic integral.

We are able to give another polynomial description of the system observing that it has
a form of D. Fairlie "elegant" integrable system (\cite{Fa}) for $n=4$ (we can skip
the unimportant coordinate $x_6$):

\begin{eqnarray}
\label{elegant}
& {\dot x_1} = x_2x_3x_4\\
& {\dot x_2} = x_1x_3x_4\\
& {\dot x_3} = x_1x_2x_4\\
& {\dot x_4} = g^2x_1x_2x_3.
\end{eqnarray}

This system admits a beautiful description as a {\it decoupled} pair
of Euler -Nahm tops after the following change of variables:
$$
\begin{array}{rcl}
 u_+ = x_3x_4 + g^2x_1x_2\\
 v_+ = x_2x_4 + g^2x_1x_3 \\
 w_+ = x_1x_4 + g^2x_3x_2\\
 u_- = x_3x_4 - g^2x_1x_2\\
 v_- = x_2x_4 - g^2x_1x_3 \\
 w_- = x_1x_4 - g^2x_3x_2 .
\end{array}
$$
In the terms of the new variables the DELL system (\ref{elegant}) is equivalent to
\begin{eqnarray}
\label{Nahm_{+}}
& {\dot u_{+}} = v_{+}w_{+}\\
& {\dot v_{+}} = w_{+}u_{+}\\
& {\dot w_{+}} = u_{+}v_{+}.
\end{eqnarray}
and to
\begin{eqnarray}
\label{Nahm_{-}}
& {\dot u_{-}} = v_{-}w_{-}\\
& {\dot v_{-}} = w_{-}u_{-}\\
& {\dot w_{-}} = u_{-}v_{-}.
\end{eqnarray}

Geometrically this change of variables means a passage to the direct ("decoupled")
product of two elliptic curves ${\cal E_{+}}\times{\cal E_{-}}$ given by the Casimirs
of the models \ref{Nahm_{+}} and \ref{Nahm_{-}}

\begin{eqnarray}
{\cal E_{+}}:
& u_{+}^2 - v_{+}^2 =k^2(E^2-1) \\
& u_{+}^2 - w_{+}^2 =(k^2 -2)(E^2 -1)
\end{eqnarray}
and
\begin{eqnarray}
{\cal E_{-}}:
& u_{-}^2 - v_{-}^2 =k^2(E^2-1) \\
& u_{-}^2 - w_{-}^2 =(k^2-2)(E^2-1).
\end{eqnarray}

This result reminds the theorem of R.Ward (\cite{ward}) that the second order
differential operator with Lamè potential $n(n+1)/2 sn^2(q|k)$ is decoupled to a
couple of the first order matrix operators $\pa + A, {\pa} - A $ if the matrix $A =
(A_1,A_2,A_3)$ satisfies to the Nahm system

$$
{\dot A_i} = \epsilon_{ijk}[A_j,A_k].
$$

\subsection{DELL as an example of Nambu - Hamilton system}

The polynomial Poisson structure describing the DELL system provides also an example
of the Nambu -Poisson structure ({\bf Subsection 1.1}).

The constructions involving quadrics gives a non-trivial example of the Nambu -
Hamilton dynamical system. Namely, the system of three quadrics in
${\CC}^5$
$$
\begin{array}{rcl}
 x_1^2 - x_2^2 =1 \\
 x_1^2 - x_3^2 =k^2 \\
 -g^2 x_1^2 + x_4^2  +x_5^2  = 1.
\end{array}
$$
admits a section by the choice of the level $x_5 = E$ and the 1-parametric
intersection of three quadrics in ${\CC}^4$
$$
\begin{array}{rcl}
 Q_1 = x_1^2 - x_2^2 =1 \\
 Q_2 = x_1^2 - x_3^2 =k^2 \\
 Q_3 = -g^2 x_1^2 + x_4^2 = 1 - E^2.
\end{array}
$$
defines the Nambu -Hamilton system (which is nothing but our old DELL system!):
$$
{{dx_i}\over{dt}} = \bigl\{Q_1,Q_2,Q_3,x_i\bigr\}.
$$

\section{Acknowledgements}
We are thankful to H. Braden and A.Gorsky for the collaboration and stimulating
discussions. Our special thanks to A. Morozov whose interest was
very inspiring to us. The numerous discussions of V.R. with A. Marshakov,
A. Mironov, N.Nekrasov  and A.Turbiner were very helpful. 
Both authors are indebted to B. Enriquez for his friendly participation
and useful discussions.

A.O. is grateful to University of Angers where the
 big part of this paper was written for invitation and very stimulating
 working atmosphere. He had acknowledged the partial support of the grants
RFBR 99-01-01669, RFBR 00-15-96579,CRDF RP1-2254 and INTAS 00-00055.
V.R. is grateful to Boris Khesin for helpful discussions  and warm
 hospitality during his stay at Toronto.
 Work of V.R. is supported partially by INTAS 99-1705, RFBR 01-01-00549 and
by the grant for support of scientific schools 00-15-96557.

\end{document}